\documentclass[12pt]{amsart}

\usepackage[OT2, T1]{fontenc}

\usepackage{amscd}
\usepackage{amsmath}
\usepackage{amssymb}
\usepackage{mathrsfs} 			
\usepackage{units}
\usepackage[all]{xy}

\usepackage{algorithm}
\usepackage{algorithmic}

\def\frk{\frak}               

\def\qq{{\frk q}}

\def\mm{{\frk m}}

\def\Phi{{\frk n}}
\def\Phi{{\frk N}}
\def\opn#1#2{\def#1{\operatorname{#2}}} 
\opn\chara{char} \opn\length{\ell} \opn\pd{pd} \opn\rk{rk}
\opn\projdim{proj\,dim} \opn\injdim{inj\,dim} \opn\rank{rank}
\opn\depth{depth} \opn\sdepth{sdepth} \opn\fdepth{fdepth}
\opn\grade{grade} \opn\height{height} \opn\embdim{emb\,dim}
\opn\codim{codim}  \opn\min{min} \opn\max{max}

\opn\Tr{Tr} \opn\bigrank{big\,rank}
\opn\superheight{superheight}\opn\lcm{lcm}
\opn\trdeg{tr\,deg}
\opn\reg{reg} \opn\lreg{lreg} \opn\ini{in} \opn\lpd{lpd}
\opn\size{size}
\opn\div{div} \opn\Div{Div} \opn\cl{cl} \opn\Cl{Cl}
\opn\Spec{Spec} \opn\Supp{Supp} \opn\supp{supp} \opn\Sing{Sing}
\opn\Ass{Ass} \opn\Min{Min}
\opn\Ann{Ann} \opn\Rad{Rad} \opn\Soc{Soc}
\opn\Im{Im} \opn\Ker{Ker} \opn\Coker{Coker} \opn\Am{Am}
\opn\Hom{Hom} \opn\Tor{Tor} \opn\Ext{Ext} \opn\End{End}
\opn\Aut{Aut} \opn\id{id}  \opn\deg{deg}

\opn\nat{nat}
\opn\pff{pf}
\opn\Pf{Pf} \opn\GL{GL} \opn\SL{SL} \opn\mod{mod} \opn\ord{ord}
\opn\Gin{Gin} \opn\Hilb{Hilb}
\opn\aff{aff} \opn\con{conv} \opn\relint{relint} \opn\st{st}
\opn\lk{lk} \opn\cn{cn} \opn\core{core} \opn\vol{vol}
\opn\link{link} \opn\star{star}
\opn\gr{gr}

\def\pot#1#2{#1[\kern-0.28ex[#2]\kern-0.28ex]}

\opn\dirlim{\underrightarrow{\lim}}
\opn\inivlim{\underleftarrow{\lim}}
\let\iso=\cong

\let\to=\rightarrow

\def\Implies{\ifmmode\Longrightarrow \else
        \unskip${}\Longrightarrow{}$\ignorespaces\fi}
\def\implies{\ifmmode\Rightarrow \else
        \unskip${}\Rightarrow{}$\ignorespaces\fi}
\def\iff{\ifmmode\Longleftrightarrow \else
        \unskip${}\Longleftrightarrow{}$\ignorespaces\fi}

\let\:=\colon

\DeclareMathOperator{\Ir}{Irr}
\DeclareMathOperator{\ch}{char}
\DeclareMathOperator{\card}{card}

\newtheorem{Theorem}{Theorem}[]
\newtheorem{Lemma}[Theorem]{Lemma}
\newtheorem{Corollary}[Theorem]{Corollary}
\newtheorem{Proposition}[Theorem]{Proposition}

\theoremstyle{definition}

\newtheorem{Remark}[Theorem]{Remark}

\newtheoremstyle{subsection-tweak}
   {11pt}
   {3pt}%
   {}
   {}%
   {\bfseries}
   {}%
   {.5em}
   {\thmnumber{\@{#1}{}\@{#2}.}%
    \thmnote{~{\bfseries#3.}}}    

\newcounter{numberingbase}

\theoremstyle{subsection-tweak}
\newtheorem{bpp}[Theorem]{}
\newtheorem{bppt}[numberingbase]{}
\newcommand{\bbpp}{\begin{bpp}}
\newcommand{\eepp}{\end{bpp}}
\newcommand{\bbppt}{\begin{bppt}}
\newcommand{\eeppt}{\end{bppt}}

\theoremstyle{theorem}

\theoremstyle{definition}

\newcommand{\val}{\mathrm{val}}		

\providecommand{\qxq}[1]{\quad\text{#1}\quad}

\newcommand{\tst}{\textstyle}

\newcommand{\sU}{\mathscr{U}}

\DeclareMathOperator{\Frac}{Frac}		
\newcommand{\wt}{\widetilde}

\let\epsilon\varepsilon
\let\phi=\varphi
\def\qed{\ifhmode\textqed\fi
      \ifmmode\ifinner\quad\qedsymbol\else\dispqed\fi\fi}
\def\textqed{\unskip\nobreak\penalty50
       \hskip2em\hbox{}\nobreak\hfil\qedsymbol
       \parfillskip=0pt \finalhyphendemerits=0}
\def\dispqed{\rlap{\qquad\qedsymbol}}

\opn\dis{dis}
\def\pnt{{\raise0.5mm\hbox{\large\bf.}}}

\opn\Lex{Lex}

\begin{document}

\title{Valuation rings of mixed characteristic as limits of complete intersection rings}

\author{ Dorin Popescu}

\dedicatory{In the memory of Richard Swan}

\address{Simion Stoilow Institute of Mathematics of the Romanian Academy,
Research unit 5, P.O. Box 1-764, Bucharest 014700, Romania,}

\address{University of Bucharest, Faculty of Mathematics and Computer Science
Str. Academiei 14, Bucharest 1, RO-010014, Romania,}

\address{ Email: {\sf dorin.m.popescu@gmail.com}}

\begin{abstract} 
Let  $V$ be 
a mixed characteristic valuation ring 
with value group $\Gamma$, valuation val, and
residue field $k$ of characteristic $p > 0$. We show that if
$\Gamma/{\bf Z} \val(p)$
has no $p$-torsion and
$V$ is Henselian, then $V$ is a filtered colimit of complete intersection ${\bf Z}_{(p)}$-algebras.


\end{abstract}

\maketitle

\section*{Introduction}

In \cite{P} we gave a different proof to a weak form of Zariski's Uniformization Theorem \cite{Z} using the following result (see \cite[Theorem 21]{P}).
\begin{Theorem}(\cite{P})\label{T0} Let $V\subset V'$ be an immediate extension of valuation rings   containing $\bf Q$.
Then $V'$ is a filtered  colimit of smooth $V$-algebras.
\end{Theorem}

 An extension of the above theorem stated in  \cite[Theorem 2]{P1} is the following theorem. 

\begin{Theorem}(\cite{P1}) \label{T0'} Let $V\subset V'$ be an extension of valuation rings containing $\bf Q$, $\Gamma\subset \Gamma'$ the value group extension of $V\subset V'$  and $\val: K'^{*}\to \Gamma'$ the valuation of $V'$. Then $V'$  is a filtered colimit of smooth $V$-algebras if and only if  the following statements hold
\begin{enumerate}

\item for each $\qq\in \Spec V$ the ideal $\qq V'$ is prime,

\item  for any prime ideals $\qq_1,\qq_2\in \Spec V$ such that $\qq_1\subset \qq_2$ and $\height(\qq_2/\qq_1)=1$  and any $x'\in \qq_2V'\setminus \qq_1'$ there exists $x\in V$ such that $\val(x')=\val(x)$,  where $\qq_1'\in \Spec V'$ is the prime ideal corresponding to the maximal ideal of $V_{\qq_1}\otimes_V V'$, that is the maximal prime ideal of $V'$ lying on $\qq_1$.
\end{enumerate}
\end{Theorem} 

If the characteristic of the residue field of $V$ is positive then  $V'$  could be not a filtered colimit of smooth $V$-algebras   shown for example in \cite[Example 3.13]{Po} (see also \cite[Remark 6.10]{Po})
 inspired by \cite[Sect 9, No 57]{O}.  This example is an algebraic extension which is not dense. If $\dim V=1$ and  the completion of $V$ is separable and transcendental over $V$ then the  immediate algebraic extension $V\subset V'$ is  dense if $V'$ is   a filtered  colimit of smooth $V$-algebras (see \cite[Theorem 2]{P''}). A {\em dense} extension  of valuation rings $V_1\subset V_2$ means that every element of $V_2$ is the limit of a fundamental sequence over $V_1$.
  Thus  $V\subset V'$ above is not a filtered  colimit of smooth $V$-algebras.
 
  If the characteristic of $V$ is positive, we  have the following result \cite[Theorem 4]{P3} (see also \cite{P6}).

\begin{Theorem}(\cite{P3})\label{T1}
Let $V\subset V'$ be an  immediate extension of valuation rings containing a field of positive characteristic and $K\subset K'$ its fraction field extension. If $K'=K(x)$ for some algebraically independent system of elements $x$ over $K$ 
then $V'$ is a filtered union of its smooth $V$-subalgebras. 
\end{Theorem}

 A form of   Zariski's Uniformization Theorem in a  case of
 positive characteristic  is a consequence of Theorem \ref{T1}.

\begin{Corollary}\cite[Corollary 5]{P3}
 Let $V$ be a  valuation ring containing its residue field $k$   with a  value group $\Gamma$ free as a $\bf Z$-module, $\val$ its valuation and $K$ its fraction field. Assume that  $K=k(x,y)$ for some  algebraically independent elements $x,y$,  $x=(x_i)_{i\in I},y=(y_j)_{j\in J}$ over $k$ such that $\val(y)$ is a basis in $\Gamma$. 
 Then $V$ is a filtered  union of its smooth $k$-subalgebras.
\end{Corollary}

An extension of this corollary is given below. 
 
\begin{Theorem}(\cite[Theorem 15]{P5}) \label{t3} Let $V\subset V'$ be an immediate extension of valuation rings of positive characteristic, $\val$ the valuation of $V'$ 
and $K\subset K'$ their fraction field extension. Assume that 
\begin{enumerate}

\item $V$ contains its residue field $k$ and $x$ is a system of elements of $V$ such that $\val(x)$ forms a $\bf Z$-linear basis of the value group of $V$.

\item $V$ is a separable dense extension of $V\cap k(x)$.

\item $K'$ is a separable dense extension of a pure transcendental field extension of $K$.
\end{enumerate}

Then $V'$ is a filtered colimit of smooth $k$-algebras.

\end{Theorem}

   In general we have  the following result  (the algebraic case is given  in \cite[Theorem 1]{P2}, see especially the arxiv version, also see \cite[Theorem 3.8]{P4}).

\begin{Theorem}(\cite[Theorem 6]{P3})\label{T2} Let  $ V'$ be an  immediate extension  of a valuation ring $V$.  Then $V'$ is a filtered
 union of its complete intersection $V$-subalgebras.
\end{Theorem}

A {\em complete intersection} $V$-algebra   is a  $V$-algebra of type $C/(P)$, where $C$ is  a polynomial $V$-algebra and $P$ is a regular system of elements of $C$ (see also \cite[Definition 19.3.1]{G}). Sometimes we call  complete intersection algebra an essentially complete intersection algebra, that is a localization of a complete intersection algebra as above. 

\begin{Remark}\label{r0} The algebraic case of Theorem \ref{T2} appeared   in a different form in  \cite[Theorem 6.2]{KT} following \cite{T1}, \cite{T2}.
\end{Remark} 
 
\begin{Remark}\label{r} Actually, Theorem \ref{T1} and  Theorem \ref{T2} are stated in \cite{P3} when $V$ contains a field. The algebraic case of  Theorem \ref{T2} given in  \cite{P2} has not this restriction. In \cite{P3} we were just interested to state Theorem  \ref{T1} when $V$ contains a field  but the proof holds also without this restriction. In fact,  the corresponding proof is easier since the characteristic of $V$ is zero. Clearly, in \cite[Lemmas 11, 13-15]{P3} it is not necessary that $V$ contains a field. This problem  appears in \cite[Proposition 16.18]{P3} but   we get new propositions in the mixed case removing the condition given when  $\ch K=p>0$. In the proof of \cite[Lemma 22]{P3} in the mixed characteristic we apply the new \cite[Proposition 16]{P3} instead \cite[Corollary 17]{P3}. Also in the proof of \cite[Theorem 23]{P3} in the mixed characteristic we apply the new \cite[Proposition 18]{P3} instead of \cite[Corollary 20]{P3}.
 \end{Remark}
 
 \begin{Theorem} Let $W\subset V$ be an immediate extension of mixed characteristic valuation rings, $p>0$ their residue field characteristic and $L\subset K$ its fraction field extension. Assume that $K$ is purely transcendental over $L$ and $W$ is absolute integrally closed  valuation ring,  Then $V$ is a filtered union of its regular ${\bf Z}_{(p)}$-subalgebras.
 \end{Theorem}
For the proof note that  $W$ is a filtered union of its regular ${\bf Z}_{(p)}$-subalgebras by \cite[Corollary 4.2.4]{AD} and so it is enough to apply Theorem
\ref{T1} (see also Remark \ref{r}).

Using Theorem \ref{T1} and Remark \ref{r}
as in \cite[Proposition 18]{P1} (see also \cite[Corollary 24]{P5}, \cite[Proposition 1]{P4}) we have the following consequence.

\begin{Corollary} \label{c5} Let $V\subset V'$ be an  immediate extension of valuation rings  and $K\subset K'$ its fraction field extension. If $V$ is Henselian and $K'=K(x)$ for some algebraically independent system of elements $x$ over $K$ then every finite system of polynomials over $V$, which has a solution in $V' $ has also one in $V$.
\end{Corollary}

The ideas of this corollary were used to show some conjectures of M. Artin in 
\cite[Theorems 1.3, 1.4]{Po0} (see also \cite{S} and \cite[Theorem 5.3.1]{I}), partially the Bass-Quillen conjecture in \cite{Po-1} (see also \cite{S}) and the Quillen conjecture in \cite{Po-2} (see also \cite{S}) using the so called the General N\'eron desingularization (see  \cite[Theorem 2.5]{Po0}, \cite[Theorem 1.1]{S} and \cite[Theorem 5.2.56]{I}).

By Theorem \ref{T2} we obtained in \cite[Theorem 9]{P5} the following result.

\begin{Theorem}(\cite{P5})\label{T3} Let $V$ be a  valuation ring containing a perfect field $F$ of positive characteristic, $k$ its residue field and   $\Gamma$ its  value group.  Then $V$ is a filtered
colimit of    complete intersection $F$-algebras  if one of the following conditions holds:
\begin{enumerate}

\item  $k\subset V$, 
\item  $V$ is Henselian. 
\end{enumerate}
\end{Theorem}

The goal of this paper is to extend the above result for the mixed characteristic valuation rings (see  Proposition \ref{p5}, Theorem \ref{t9}). 

\begin{Theorem}\label{T4} Let $V$ be a   mixed characteristic valuation ring. Then $V$ is a filtered colimit of 
 complete intersection ${\bf Z}_{(p)}$-algebras 
if either of the following conditions hold:
\begin{enumerate}
\item  $V$ has a  DVR subring $(R,\pi)$ with the same residue field such that $\pi V$ is a prime ideal. 
\item  $V$ is Henselian and  $\Gamma/{\bf Z}\val(p)$ has no $p$-torsion, where $\Gamma$ is the value group of $V$, $\val$ its valuation and $p$ is the characteristic of its residue field. 
\end{enumerate} 
\end{Theorem}

An abelian group $(G,+)$ has {\em no $p$-torsion} if no nonzero element of it could be  killed  by powers of $p$. 
\vskip 0.3 cm

As a consequence some Andr\'e-Quillen homology and cohomology (see \cite{An}) of a valuation ring are zero, if  it is  subject to one of the above theorems  (see Corollary \ref{hom}).

We owe thanks to an anonymous referee who almost hinted us Theorem 9 
and had some useful comments especially on  the proofs of Proposition 20 and  Theorem 22.

\section{Valuation rings of mixed characteristic with finitely generated value group}

A field extension $K\subset K'$ is {\em separably generated} if $K'$ is an algebraic separable extension of a pure transcendental extension of $K$. Thus a separable finite type field extension is separably generated.

\begin{Lemma} \label{l} Let $V$ be a Henselian mixed characteristic valuation ring, $k$ its residue field and $p=\ch(k)$. Assume that $k'\subset k$ is a subfield  separably generated over ${\bf F}_p$ the finite field with $p$ elements. Then there exists a DVR subring $R\subset V$ such that
\begin{enumerate}
\item $pR$ is the maximal ideal of $R$,
\item $R$ is a filtered union of smooth ${\bf Z}_{(p)}$-subalgebras. 
\item $R\subset V$ is an extension of valuation rings with $k'\subset k$ the  residue field extension.
\item if $ k'$ is finitely generated over ${\bf F}_p$ then $R$ is  essentially of finite type over ${\bf Z}_{(p)}$. 
\end{enumerate}
\end{Lemma}
{\renewcommand{\qedsymbol}{}
\begin{proof}  By hypothesis there exists a system of elements $x$ of $V$ inducing a separable transcendental basis of $k'$ over ${\bf F}_p$,
that is $k'$ is algebraic separable over ${\bf F}_p(x)$ . Then $V_1=({\bf Z}_{(p)}[x])_{p{\bf Z}_{(p)}[x]}$ is a DVR, subring of $V$, $pV_1$ is its maximal ideal and the residue field extension of $V_1\subset V$ is algebraic separable. Let ${\bar y}\in k'$ which is not in the residue field $k_1$ of $V_1$ and ${\bar f}= \Ir({\bar y},k_1)\in k_1[Y]$, that is the minimal polynomial with root $\bar y$ from $k_1[Y]$. Let $f\in V_1[Y]$ be a monic polynomial lifting $\bar f$. As $V$ is Henselian we may lift $\bar y$ to a solution $y$ of $f$ in $V$ and $V_2=(V_1[Y]/(f))_{pV_1[Y]}$ is a DVR, which can be considered a subring of $V$ via the injective map $V_2\to V$ extending $V_1\subset V$ by $Y\to y$. Clearly,  $\bar y$ is contained in the residue field of $V_2$. Using this trick by transfinite induction or by Zorn's Lemma  we find such $R$.  

If $ k'$ is finitely generated over ${\bf F}_p$ then above 
$x$ is a finite set and we may take just one $\bar y$ and so $R$ is essentially of finite type over ${\bf Z}_{(p)}$.

 Note that $R$ is a  filtered  union of its smooth ${\bf Z}_{(p)}$-subalgebras by N\'eron desingularization  \cite{N} (see also \cite{KAP}).\hfill$\square$
\end{proof}
}

\begin{Proposition} \label{p} Let $V$ be a Henselian mixed characteristic valuation ring, $k$ its residue field, $\ch k=p$, $\val$ its valuation and $\Gamma$ its value group. Assume that $\Gamma$ is finitely generated, $\height(pV)=1$,
 $\Gamma/{\bf Z}\val(p)$ is torsion free  and  $k$ is separably generated over ${\bf F}_p$. Then  $V$ is a filtered  union of its complete intersection subalgebras over a DVR  $R$ with its maximal  
 ideal $pR$ and its  residue field $k$. Moreover, $V$ is a filtered  union of its complete intersection subalgebras over 
 ${\bf Z}_{(p)}$.
\end{Proposition}
\begin{proof}  By Lemma \ref{l} there exists a DVR subring $R\subset V$ such that
 $pR$ is the maximal ideal of $R$, 
 the extension $R\subset V$ has the trivial residue field extension and $R$ is  a  filtered  union of its smooth ${\bf Z}_{(p)}$-subalgebras $R_{\alpha}$. 
 Clearly, $\Gamma/{\bf Z}\val(p)$ is free because it is finitely generated and has no torsion.
 Choose a $\bf Z$ basis induced by $\val(x)$, for some  $x$ from $V$. Thus $x$ is algebraically independent over $R$ by \cite[Theorem 1 in VI (10.3)]{Bou}. Let $T$ be the fraction field of $R$ and set $W=V\cap T(x)$.  By \cite[Lemma 26 (2)]{P},  $ W$ is a   filtered  union of its complete intersection $R$-subalgebras $W_{\beta}$, even regular local rings and $pW$ is a prime ideal.  Note that $\Gamma$ is the value group of $W$, that is the extension $W\subset V$ is immediate and so $pV$ is also prime.  Using Theorem \ref{T2} we see that $V$ is a   filtered  union of its complete intersection $W$-subalgebras $V_{\gamma}$ since  $W\subset V$ is immediate. 
 
 For each $\beta$ there exist $\alpha $ and a complete intersection $R_{\alpha}$-subalgebra $W_{\beta \alpha}$ of $R$ such that
 $W_{\beta}\iso R\otimes_{R_{\alpha}} W_{\beta \alpha}$. But $W_{\beta \alpha}$  is a complete intersection ${\bf Z}_{(p)}$-algebra by \cite[Lemma 6]{P2}. Thus $W$ is a filtered union of complete intersection ${\bf Z}_{(p)}$-subalgebras $W_{\beta \alpha}$ of $W$.
 
 Similarly, for each $\gamma$ there exist $\beta, \alpha$ and a
 complete intersection $W_{\beta \alpha}$-subalgebra $V_{\gamma \beta \alpha}$ of $V$ such that $V_{\gamma}\iso W\otimes_{W_{\beta \alpha}} V_{\gamma \beta \alpha}$ and using again  \cite[Lemma 6]{P2} we see that $V_{\gamma \beta \alpha}$ is a  complete intersection ${\bf Z}_{(p)}$-algebra. 
 
 Thus $V$ is  a filtered union of complete intersection  ${\bf Z}_{(p)}$-subalgebras  of $V$ using \cite[Lemma 1.5]{S}.
\hfill\ \end{proof}

\begin{Lemma}\label{l'} Let $R\subset V$ be an extension of valuation rings
with the same residue field $k\supset {\bf F}_p$, $T\subset K$ its fraction field extension,  $\val$ the valuation of $V$ and $\Gamma$ the value group of $V$. Assume that $V$ is Henselian, $R$ is a DVR with $\pi$ its local parameter,
$\Gamma$ is finitely generated and $\Gamma/{\bf Z}\val(\pi)$ has no $p$-torsion. Then there exists a DVR $A\subset V$ containing   $R$ with $\pi_1$ its local parameter such that the inclusions $R\subset A\subset  V$ are extensions of valuation rings, the extension  $R\subset A$ is a filtered union of complete intersection $R$-subalgebras   and $\Gamma/{\bf Z}\val(\pi_1)$ is torsion free. 
\end{Lemma}
\begin{proof} We may assume that $\Gamma/{\bf Z}\val(\pi)$ is not torsion free, otherwise take $A=R$. Let $\nu\in \Gamma$ which is not in ${\bf Z}\val(\pi)$ but $t\nu\in {\bf Z}\val(\pi)$ for some $t\in {\bf N}\setminus p{\bf N}$, $t>1$. Thus $t\nu=n\val(\pi)$ for some $n\in {\bf Z}$. We may assume that gcd$(t,n)=1$. Then we have $1=at+bn$ for some $a,b\in {\bf Z}$ and we get $t(b\nu+a\val(\pi))=\val(\pi)$. Consequently for $\nu'=b\nu+a\val(\pi)$ we have $t\nu'=\val(\pi)$ and so $d\pi= y^t$ for some $y\in V$
and an unit $d$ of $V$. 
As the residue field extension of $R\subset V$ is trivial we have $d=cd'$ for some units $c\in R$ and $d'\in V$ with $d'\equiv 1$ modulo $\mm$ the maximal ideal of $V$.  The equation $Y^t-d'=0$ has the solution $1$ modulo $\mm$, which can be lifted to a solution $d''$ in $V$ by the Henselian property (note that $t\not \in p{\bf Z}$) . Set $\pi_1=y/d''$ and we get $\pi_1^t=c\pi$.  Note that the DVR $A=V\cap T(\pi_1)\cong (R[Y]/(Y^t-c\pi))_{(Y)}$ is a   complete intersection, essentially of finite type $R$-algebra. 

Changing  from $R$ to $A$ and $\pi$ to $\pi_1$ the torsion part of $\Gamma/{\bf Z}\val(\pi_1)$ becomes smaller and using this trick step by step, $\Gamma$ being finitely generated, we arrive to a DVR $A'$ with its local parameter $\pi'$ such that $\Gamma/{\bf Z}\val(\pi')$ is torsion free.
\hfill\ \end{proof}

\begin{Theorem}\label{t8} Let $V$ be a Henselian mixed characteristic valuation ring, $k$ its residue field, $p=\ch k$, $\val$ its valuation and $\Gamma$ its value group. Assume that $k$ is separably generated over ${\bf F}_p$,  $\Gamma$ is finitely generated and $\Gamma/{\bf Z}\val(p)$ has no $p$-torsion.
  Then  $V$ is a filtered  union  of its complete intersection subalgebras over a DVR  $A\subset V$ with  its  residue field $k$ and its local parameter $\pi$ such that
 $\Gamma/{\bf Z}\val(\pi)$ has no torsion and $A\subset V$ is an extension of valuation rings. Moreover, $V$ is a  filtered  union  of its complete intersection subalgebras over ${\bf Z}_{(p)}$.
\end{Theorem}
{\renewcommand{\qedsymbol}{}
\begin{proof} Let $R$ be the DVR given by Lemma \ref{l}.
By Lemma \ref{l'} there exists a DVR subring $A\subset V$ with $\pi$ its local parameter and such that the inclusion $A\subset  V$ is an extension of valuation rings, $\Gamma/{\bf Z}\val(\pi)$ is torsion free and $A$ is  a filtered   union  of its complete intersection subalgebras over ${\bf Z}_{(p)}$.

If $\height(\pi V)=1$ we can do as in the proof of Proposition \ref{p}. Assume that $ \height(\pi V)>1$. Let $\Gamma_1$ be the value group of the valuation ring $V_{\pi}$. It is free because it is finitely generated and let $z$ be a system of elements of $V$ such that $\val(z)$ is a basis of $\Gamma_1$. Also choose a $\bf Z$-basis given by $\val(x)$, $x\in V^e$ in $\Gamma/{\bf Z}\val(\pi)$. Thus $x,z$ is algebraically independent over $A$ by \cite[Theorem 1 in VI (10.3)]{Bou}  because $\val(z),\val(x)$ are linearly independent over $\bf Z$. Let $T$ be the fraction field of $A$ and set $W=V\cap T(x,z)$.  By \cite[Lemma 26 (2)]{P},  $ W$ is a   filtered  union of its complete intersection $A$-subalgebras, even regular local rings and $\pi W$ is a prime ideal.  Note that $\Gamma$ is the value group of $W$, that is the extension $W\subset V$ is immediate and so $\pi V$ is also prime.  Using Theorem \ref{T2} we see that $V$ is a   filtered 
 union  of its complete intersection $W$-subalgebras, since  $W\subset V$ is immediate. As in the proof of Proposition \ref{p} we see that $V$ is a   filtered 
 union  of its complete intersection $W$-subalgebras of $V$ using 
\cite[Lemma 6]{P2} and \cite[Lemma 1.5]{S}.  
\hfill$\square$\end{proof}
}

The following proposition is necessary in the next section.

\begin{Proposition} \label{p4} A mixed characteristic DVR
 $V$ with its residue field  separably generated over ${\bf F}_p$ is a filtered colimit of some complete intersection algebras over ${\bf Z}_{(p)}$.
\end{Proposition}

\begin{proof} Let $x$ be a system of elements of $V$ lifting a transcendental basis $\bar x$ of the residue field $k$ of $V$ over ${\bf F}_p$ such that ${\bf F}_p({\bar x})\subset k$ is algebraic separable. Then $W={\bf Z}[x]_{p{\bf Z}[x]} \subset V$
is a DVR subring and the residue field extension of $W\subset V$ is algebraic separable. If $W\subset V$ is unramified then apply N\'eron Desingularization. Otherwise, as in 
\cite[Lemma 2.2]{Po1} there exists a DVR subring $B\subset V$ containing $W$ such that for a local parameter $x$ of $V$ a power 
$b=x^s$, $s>1$ is a local parameter of $B$, $C=B[x]_{(x)}\iso
(B[X]/(X^s-b))_{(X)}$  is a DVR subring of $V$  and the extension 
$C\subset V$ is unramified. Then $C$ is a complete intersection $B$-algebra and $V$ is a filtered union of its complete intersection $C$-subalgebras by N\'eron desingularization. As in Proposition \ref{p} we see that $V$ is a filtered union of its complete intersection $B$-subalgebras.  Repeat the procedure with $B$ instead $V$. Note that we will get a  power $b'=b^r$, $r<s$ as a  local parameter in $B'$ the new $B$. So in  finite  steps  the procedure ends. 
\hfill\ \end{proof}

Also we need later the following proposition.

\begin{Proposition}(\cite[Proposition 22]{P}) \label{Pr}   Let $V\subset V'$ be an extension of  valuation rings. Suppose that

\begin{enumerate}

\item $V$ is a DVR extending ${\bf Z}_{(p)}$ with $\pi$ its local parameter, and $p$ a prime number.

\item $\pi V'$ is the maximal ideal of $V'$,

\item the residue field extension of $V\subset V'$ is separable.
\end{enumerate}

Then $ V'$ is a filtered colimit of some smooth $V$-algebras.
\end{Proposition}

\vskip 0.5 cm

\section{Valuation rings of mixed characteristic with general value group}

We recall the following result from \cite{P}  obtained using methods from model theory. A {\em cross-section} of a valuation ring $V$ with value group $\Gamma$ is a section $s:\Gamma\to K^*$ of its valuation $\val:K^*\to \Gamma$ in the category of abelian
groups, following \cite[p. 293]{P}.

\begin{Theorem}(\cite[Variant A 11]{P}) \label{bt} 
For every faithfully flat map $V \subset V'$ of valuation rings  with value groups $\Gamma \subset \Gamma'$ such that $\Gamma'/\Gamma$ is torsion free, there is a countable sequence of ultrafilters ${\sU}_1, {\sU}_2, \dots$ on some respective sets $U_1, U_2, \dots$ for which the valuation rings $\{V_n\}_{n \geq 0}$ and $\{V'_n\}_{n \geq 0}$ defined inductively by $V_0 := V$ and $V'_0 := V'$ with $V_{n + 1} := \prod_{{\sU}_{n + 1}} V_{n}$ and $V'_{n + 1} := \prod_{{\sU}_{n + 1}} V'_{n}$ are such that the valuation ring 
\[
\tst  \wt{V}' = \varinjlim_{n \geq 0} V'_n \qxq{with} \wt{K}' := \Frac(\wt{V}') \qxq{has a cross-section}  \wt{s} :\wt{\Gamma}' \to \wt{K}'^*
\]
whose restriction to the value group $\wt{\Gamma}$ of $\wt{V} := \varinjlim_{n \geq 0} V_n$ lands in $\wt{K} := \Frac(\wt{V})$. 
\end{Theorem}

The next proposition is similar to Proposition \ref{p} when $V$ is not Henselian and $\Gamma$ is not necessarily finitely generated. The proof goes as in Theorem \ref{t8} because now we have by assumption the necessary  DVR $A$,  which is a filtered colimit of some complete intersection algebras over ${\bf Z}_{(p)}$ by Proposition
\ref{p4}.

\begin{Proposition} \label{p5}  Let $V$ be a  mixed characteristic valuation ring, $k$ its residue field, $p=\ch(k)$,  $\val$ its valuation and $\Gamma$ its value group. Assume that there exists a mixed characteristic DVR subring $A\subset V$ of residue field $k$ with $\pi$ a local parameter such that $A\subset V$ is an extension, $\pi V$ is a prime ideal and $\Gamma /{ \bf Z} \val(\pi)$ has no torsion. Then $V$ is a  filtered colimit of some complete intersection $A$-algebras  and in particular $V$ is a filtered colimit of some complete intersection algebras over ${\bf Z}_{(p)}$.
\end{Proposition}

\begin{proof} Apply Theorem \ref{bt} for the faithfully flat extension $A\subset V$.
Let ${\sU}_1, {\sU}_2, \dots$ be a countable set of ultrafilters on some respective sets $U_1, U_2, \dots$ for which the valuation rings
 $(V_n)_{n \geq 0}$ defined inductively by $V_0 := V$, $V_{n + 1} := \prod_{{\sU}_{n + 1}} V_{n}$ for $n \geq 1$ are such that for the valuation ring $ \wt{V} := \varinjlim_{n \geq 0} V_n$ there exists a cross-section  $\wt{s}$ of $\wt{V}$. Note that the ultrapower of a Henselian ring is always a Henselian ring, because the property of being Henselian is expressible by a first-order sentence  and it is  preserved under taking ultraproducts and ultrapowers. Thus $V_n$ and $ \wt{V}$ are Henselian.
  Similarly, we define some rings inductively by
  $A_0 := A$, $A_{n + 1} := \prod_{{\sU}_{n + 1}} A_{n}$ for $n \geq 1$ and let $ \wt{A}$  be the union of $(A_n)$. 
The residue field extension of the extension  $A_n\subset A_{n+1}$
is  separable and so the residue  field extension of the extension  $A\subset A_n$  is  separable too.  Also note that $\pi A_n$, $n\in {\bf N}$ are maximal ideals because ultrapowers of fields as $A/(\pi)$ are fields (for details concerning ultrapowers see for example \cite[page 301]{P}).  Then $\pi \wt{A}$ is the maximal ideal of $\wt{A}$ and the value group of $\wt{A}$ is the ultrapower of $\bf Z$.  We may assume by Theorem \ref{bt} that   $\wt{s}$ induces by restriction a cross-section of $\wt{A}$.

    By Proposition \ref{Pr} we see that 
    $A_n$ is a filtered colimit of some smooth $A$-algebras and so
$\wt{A}$ is a filtered  colimit of some smooth $A$-algebras. Note that $\wt{A}$ and $\wt{V}$ have the same residue field.

 Let $T$  be the fraction field of $\wt{A}$,   $\wt{\Gamma}$ the value group of $\wt{V} $ and $\Gamma'\subset \wt{\Gamma}$ a finitely generated subgroup containing $\val(\pi)$. Then $\Gamma'/{\bf Z}\val(\pi) \subset \wt{\Gamma}/{\bf Z}\val(\pi)$ has no torsion (thus free) because an ultrapower of a torsionless group has no torsion. So we may choose  a finite set $\gamma'$ in $\Gamma'$  inducing a basis in $\Gamma'/{\bf Z}\val(\pi)$. We recall that $\wt{s}$ induces by restriction a cross-section of $\wt{A}$.
  Then $\wt{s}(\gamma')$ is algebraically independent over $T$ by \cite[Theorem 1 in VI (10.3)]{Bou}. Thus $W_{\Gamma'}=V\cap T(\wt{s}(\gamma'))$ is a union of localization of complete intersection $\wt{A}$-subalgebras of $\wt{V}$ by \cite[Lemma 26 (2)]{P}. Note that $W=V\cap T(\wt{s}(\wt{\Gamma}))$ is a filtered union of  such $W_{\Gamma'}$ and so of complete intersection $\wt{A}$-subalgebras  of $\wt{V}$. As in the proof of Proposition \ref{p} we see that $W$ is a filtered 
  colimit of $A$-subalgebras of $\wt{V}$ using \cite[Lemma 6]{P2}.

 Note that  the extension $W\subset \wt{V}$ is immediate. 
Using Theorem \ref{T2} and Remark \ref{r} we see that  $\wt{V}$ is a  filtered union of some complete intersection algebras over $W$ and so over 
$\wt{A}$, even over  ${\bf Z}_{(p)}$ by  Proposition \ref{p4}.

Let $E$ be a finitely generated ${\bf Z}_{(p)}$-algebra and $w:E\to V$ a morphism. Then the composite map $E\to V\to \wt{V}$ factors through a complete intersection ${\bf Z}_{(p)}$-algebra $D$. 
Now we proceed as in \cite[Theorem 35]{P} (there we had the case of  a smooth algebra $D$ and now we  deal with  several polynomials $h$). So $w$ 
 factors through $D$ too    because all finite systems of polynomial equations, which have a solution in $\wt V$, must have one in $V$. This is enough  by \cite[Lemma 1.5]{S}.
\hfill\ \end{proof}

We need \cite[Proposition A.6]{P},  which is obtained using \cite[Theorem 6.1.4]{CK} and says in particular the following:

\begin{Proposition}\label{kes} Let $V$ be a valuation ring with value group $\Gamma$. Then there exists an ultrafilter ${\sU}$ on a set $U$ such that any system of polynomial equations
$(g_i((X_j)_{j\in J})_{i\in I}$ with $\card(I)\leq \card (U)$ in variables $(X_j)_{j\in J}$ with coefficients in  the ultrapower ${\tilde V}=\Pi_{{\sU}}V$ has a solution in  ${\tilde V}$ if and only if all its 
finite subsystems have.
\end{Proposition}

\begin{Theorem}\label{t9}  Let $V$ be a Henselian mixed characteristic valuation ring, $k$ its residue field, $p=\ch(k)$, $\val$ its valuation and $\Gamma$ its value group. Assume that   $\Gamma /{ \bf Z} \val(p)$ has no $p$-torsion. Then  $V$ is a filtered colimit of some complete intersection algebras over ${\bf Z}_{(p)}$.
\end{Theorem}
\begin{proof}
Let $\Gamma'\subset \Gamma$  a finitely generated subgroup and $k'\subset k$ a finitely generated subfield. By Lemma \ref{l} we find a DVR subring $R_{k'}\subset V$ with $k'$ its residue field, $pR_{k'}$ its maximal ideal, which is a filtered union of smooth ${\bf Z}_{(p)}$-algebras. Note that $R_{k'}$ is essentially of finite type over $\bf Z$, because $k'$ is finitely generated. Using Lemma \ref{l'} (in its proof  it is important that $V$ is Henselian)
 there exists a DVR subring $A_{k',\Gamma'}\subset V$ containing  $R_{k'}$
with a local parameter $\pi_{\Gamma'}$ such that $\Gamma'/\Gamma'\cap {\bf Z}\val(\pi_{\Gamma'})$ is torsion free and $A_{k',\Gamma'}$ is a filtered colimit of complete intersection ${\bf Z}_{(p)}$-algebras. Moreover $A_{k',\Gamma'}$ is essentially of finite type  (even essentially of finite presentation by \cite[Theorem 4]{Na}) over $\bf Z$.

Actually, $A_{k',\Gamma'}$  is essentially generated by some algebraically independent elements $x'$ over $\bf Q$ followed by some integral elements $(x''_i)_i$ over ${\bf Z}[x']$ (see the proof of Lemma \ref{l}).  We have
  $A_{k',\Gamma'}\iso{\bf Z}[x',x'']_{\mm\cap {\bf Z}[x',x'']}$, $\mm$ being the maximal ideal of $V$.
  Certainly we can choose  the $x''$ include also $\pi_{\Gamma'}$.

Choose a finite system $G_{k',\Gamma'}$ of monic polynomials in 
$X''_i$
 which yield equations of
integral dependence for the $x''_i$. 
The system $G_{k',\Gamma'}$ generates the kernel of the map
${\bf Z}[X',X'']\to V$ ,
$(X',X'')\mapsto (x', x'')$.
Actually $G_{k',\Gamma'}$  consists in some monic polynomials in $X''_i$ corresponding to $x''_i$ over ${\bf Z}[X']$. Any solution ${\bar x}', \bar{x}''$ of 
$G_{k',\Gamma'}$ in $V$ defines a DVR subring ${\bf Z}[{\bar x}',{\bar x}'']_{\mm\cap {\bf Z}[\bar{x}',\bar{x}'']}$ of $V$ isomorphic with $A_{k',\Gamma'}$.

 Let $\mathcal{E}$ be the set of all pairs $(k',\Gamma')$
with $k'\subset k$ a finitely generated subfield and $\Gamma'\subset \Gamma$ a finitely generated subgroup. For some other  $(k'',\Gamma'')\in \mathcal{E}$ with $k'\subset k''$, $\Gamma'\subset \Gamma''$ maybe $A_{k',\Gamma'}\not \subset A_{k'',\Gamma''}$. It is also possible that $A_{k',\Gamma'}$ is not contained in the Henselization $C_{k'',\Gamma''}$ of $A_{k'',\Gamma''}$. 

We claim there exists a DVR isomorphic to $A_{k',\Gamma'}$ contained in $C_{k'',\Gamma''}$. The polynomials
$G_{k',\Gamma'}$ have a solution in $C_{k'',\Gamma''}$ modulo its maximal ideal because they have a solution
in $k'$, which is contained in the residue field 
$k''$ of $C_{k'',\Gamma''}$. Because $C_{k'',\Gamma''}$ is Henselian,
there exists a solution of the $G_{k',\Gamma'}$ in $C_{k'',\Gamma''}$. By the previous paragraph, we can
therefore find a DVR isomorphic to  $A_{k',\Gamma'}$ contained in $C_{k'',\Gamma''}$.
Similar to the case of $A_{k',\Gamma'}$ we see that $C_{k',\Gamma'}$ is generated by some algebraically independent elements $y'$ over $\bf Q$ (they are present already in $A_{k',\Gamma'}$) followed by some integral elements $(y''_i)_i$ over ${\bf Z}[y']$ because $C_{k',\Gamma'}$ is defined by some etale neighborhoods of $A_{k',\Gamma'}$, which are essentially finite over 
$A_{k',\Gamma'}$ (see \cite[Theorem 2.5]{S}). Note that the extension $A_{k',\Gamma'}\subset C_{k',\Gamma'}$ is immediate.

Choose a  system $G'_{k',\Gamma'}$ (not finite) of monic polynomials in 
$Y''_i$
 which yield equations of
integral dependence for the $y''_i$. 
The system $G'_{k',\Gamma'}$ generates the kernel of the map
${\bf Z}[Y',Y'']\to V$ ,
$(Y',Y'')\mapsto (y', y'')$.
 A solution of $G'_{k',\Gamma'}$ in $V$ defines a Henselian   DVR subring of $V$
isomorphic with $C_{k',\Gamma'}$.

Let   $(k',\Gamma'),(k_1',\Gamma_1')\in \mathcal{E}$ with $k'\subset k_1'$, $\Gamma'\subset \Gamma_1'$. We claim there exists a DVR isomorphic to $C_{k',\Gamma'}$ contained in $C_{k'',\Gamma''}$. The polynomials
$G'_{k',\Gamma'}$ have a solution in $C_{k'',\Gamma''}$ modulo its maximal ideal because they have a solution
in $k'$, which is contained in the residue field 
$k'_1$ of $C_{k'_1\Gamma'_1}$. Because $C_{k'_1,\Gamma'_1}$ is Henselian,
there exists a solution of the $G'_{k',\Gamma'}$ in $C_{k'_1,\Gamma'_1}$. Therefore we can
 find a DVR isomorphic to  $C_{k',\Gamma'}$ contained in $C_{k'',\Gamma''}$.

For each $y_{k',\Gamma',i} \in  \{y', y''\}$, choose a unit
$u_{k',k'_1 ,\Gamma',\Gamma'_1,i}\in C_{k'_1,\Gamma'_1}$ such that\\
$u_{k',k'_1 ,\Gamma',\Gamma'_1,i}y_{k',\Gamma',i} 
\in {\bf Z}[y', y'']$. Then, there exist a polynomial 
$H_{k',k'_1,\Gamma',\Gamma'_1,i}$ over ${\bf Z}$  in variables  $Y_{k'_1 ,\Gamma'_1}$
corresponding
to the generator $y_{k'_1,\Gamma'_1}$
such that
$$u_{k',k'_1 ,\Gamma',\Gamma'_1,i}y_{k',\Gamma',i}=H_{k',k'_1,\Gamma',\Gamma'_1,i}(y_{k'_1 ,\Gamma'_1}).$$

 Let $F_{k',k_1',\Gamma',\Gamma_1',i}$ be the system of polynomials 
 $$U_{k',k_1',\Gamma',\Gamma_1',i} Y_{k',\Gamma',i}-H_{k',k_1',\Gamma',\Gamma_1',i}(Y_{k_1',\Gamma_1'}),$$ 
  and $U_{k',k_1',\Gamma',\Gamma_1',i}U'_{k',k_1',\Gamma',\Gamma_1',i}-1$,
  in some variables $ Y_{k',\Gamma',i}$, $ Y_{k_1',\Gamma_1',j}$, $U_{k',k_1',\Gamma',\Gamma_1',i}$,\\
  $U'_{k',k_1',\Gamma',\Gamma_1',i}$. Here the variables $ Y_{k',\Gamma',i}$ are in fact some variables corresponding to $y',y''$ as it was used above and the variables $U_{k',k_1',\Gamma',\Gamma_1',i}$ correspond to $u_{k',k'_1 ,\Gamma',\Gamma'_1,i}$.
  
  A solution of the  system of polynomials 
$G'_{k',\Gamma'}$, $G'_{k_1',\Gamma_1'}$, $(F_{k',k_1',\Gamma',\Gamma_1',i})_i$ in $V$  defines an extension of Henselian DVR's of type  $C_{k',\Gamma'}\subset C_{k_1',\Gamma_1'}$. A solution of all $G'_{k',\Gamma'}$,  $(F_{k',k_1',\Gamma',\Gamma_1',i})_i$ in $V$  defines a filtered set by inclusion of Henselian DVR's corresponding to $(k',\Gamma'), (k_1',\Gamma_1')\in \mathcal{E}$.

Let $T_{k',\Gamma'} $ be the fraction field of  $C_{k',\Gamma'}$
for some $(k',\Gamma')\in \mathcal{E}$.  Let  $\Gamma_1$ be the valuation group of the localization $V[1/p]$ of $V$ with respect to the powers of $p$. By \cite[Ch. VI, (4,3)]{Bou} there exists a surjective canonical map $\rho:\Gamma\to \Gamma_1$. As $\rho(\Gamma') $ is finitely generated and so free,  choose some elements $\gamma'$ in $\Gamma'$ such that $\rho(\gamma')$ is a basis in $\rho(\Gamma')$. Also choose some elements $\gamma''$ in $\Gamma'$ which lift a basis of $\Gamma'/\Gamma' \cap {\bf Z} \val(\pi_{\Gamma'})$. Then $\gamma', \gamma''$ is a $\bf Z$- linearly independent set in $\Gamma'$.  
Choose a system of elements $z'_{\Gamma'}$ of $V$ such that $\gamma'=\val(z'_{\Gamma'})$ and $z''_{\Gamma'}$ a system of elements of $V$ such that $\gamma''=\val(z''_{\Gamma'})$. Then 
$z'_{\Gamma'}, z''_{\Gamma'}$ are algebraically independent over 
$T_{k',\Gamma'}$ (see 
\cite[Theorem 1 in VI (10.3)]{Bou}). Note that 
$W_{k',\Gamma'}=V\cap T_{k',\Gamma'}(z'_{\Gamma'},z''_{\Gamma'})$
is a valuation ring with $\pi_{\Gamma'}W_{k',\Gamma'}$ a prime ideal, its residue field $k'$ and its value group $\Gamma'$. Moreover  $W_{k',\Gamma'}$ is a filtered colimit of complete intersection  ${\bf Z}_{(p)}$-algebras (see the proof of Theorem  \ref{t8} and \cite[Lemma 26 (2)]{P}).

 Let $W_{k',\Gamma'}^h\subset V$ be the Henselization of $W_{k',\Gamma'}$. Note that $W_{k',\Gamma'}^h$ is still a filtered colimit of complete intersection  ${\bf Z}_{(p)}$-algebras. Let  $L_{k',\Gamma'}$ be a system of polynomials over more variables $Y$, which include also
  $G'_{k',\Gamma'}$, such that $W_{k',\Gamma'}^h\iso {\bf Z}[Y]/(L_{k',\Gamma'})$.  For some other  $(k'',\Gamma'')\in \mathcal{E}$ with $k'\subset k''$, $\Gamma'\subset \Gamma''$ we may find in  $W_{k'',\Gamma''}^h$ a DVR 
 of type $A_{k',\Gamma'}$ as above and even a valuation ring of type $W_{k',\Gamma'}^h$. Certainly, the inclusion $W_{k',\Gamma'}^h \subset W_{k'',\Gamma''}^h$ means some other equations $S_{k',k'',\Gamma',\Gamma''}$ in more variables $Y,Y', U,U'$ extending the previous $F_{k',k'',\Gamma',\Gamma''}$.

  We apply Proposition \ref{kes}. There exists an ultrafilter $\mathcal{P}_1$ on a set 
  $P_1$ with $\card(P_1)$ greater than the cardinal of  all the systems $L_{k',\Gamma'}$, $S_{k',\Gamma'}$ (in particular also
  $G'_{k',\Gamma'}$,  $(F_{k',k_1',\Gamma',\Gamma_1',i})_i$) such that $L=(L_{k',\Gamma'}) $ and $S=(S_{k',k_1',\Gamma',\Gamma_1'})$ have  a solution in the ultraproduct
    $V_1=\Pi_{\mathcal{P}_1}  V$ because each finite subsystem of them has a solution in $V_1$ (even in $V$). Indeed, for a finite set of elements  $({\bar k}'_j,{\bar\Gamma}'_j)$ of $\mathcal{E}$     
    we take a $({\bar k},\bar{\Gamma})$ of  $\mathcal{E}$ such that ${\bar k}$ contains all ${\bar k}'_j$ and $\bar{\Gamma}$ contains all $\bar{\Gamma}'_j$. 
  Then the valuation ring given by the solution of 
 $L_{\bar{k},\bar{\Gamma}} $ contains a solution of all $L_{\bar{k}_j,\bar{\Gamma}_j}$, $S_{\bar{k}_j,\bar{k}_{j'},\bar{\Gamma}_j,\bar{\Gamma}_{j'}}$.
 
 Fix such a solution of $L $, $S$  in $V_1$.   $V_1$ is  Henselian as in the proof of Proposition \ref{p5}. So there exists a filtered set by inclusion of Henselian valuation rings $W_{k',\Gamma'}^h$ in particular also of Henselian DVR's subrings $C_{k',\Gamma'}$ of $V_1$, which are filtered colimits of  complete intersection ${\bf Z}_{(p)}$-algebras (see example Proposition \ref{p4}). The union $W_1\subset V_1$ of $(W_{k',\Gamma'}^h)$ is a Henselian valuation ring  with its residue field $k$ and with residue field $\Gamma$. Also the union $C_1$ of $C_{k',\Gamma'}$   is  a valuation ring  with its maximal ideal  the radical of $pC_1$.

   Repeating this procedure with $V_1$ instead $V$  we find a set $P_2$ and an ultrafilter $\mathcal{P}_2$ such that $V_2=\Pi_{\mathcal{P}_2} V_1$ contains a   valuation ring $W_2$  with the  residue field $k_1 $ of $V_1$ and with its value group $\Gamma_1$  the value group of $W_1$, which is a filtered colimit of  complete intersection ${\bf Z}_{(p)}$-algebras. Repeating again this procedure we find
some sets $(P_n)_n$ and some   ultrafilters $(\mathcal{P}_n)_n$ on them and define $V_{n+1}=\Pi_{\mathcal{P}_{n+1}} V_n$ and  $\wt V=\varinjlim_{n \geq 0} V_n$. In this way we obtain a filtered set  ordered  by inclusion by pairs $W_{n+1}\subset V_{n+1}$ of   valuation rings 
 ordered by $(W_{n+1}\subset V_{n+1})\leq (W_{n+2}\subset V_{n+2})$ if $W_{n+1}\subset W_{n+2}$, $V_{n+1}\subset V_{n+2}$.
Note that the  residue field of $W_{n+1}$ is  $k_n $ the residue field  of $V_n$ and the  value group of $W_{n+1}$ is the value group  of  $V_n$. Moreover, they are filtered colimits of  complete intersection
 ${\bf Z}_{(p)}$-algebras.

 So the union $W_{\infty}$ of $W_n$ is a filtered colimit of complete intersection ${\bf Z}_{(p)}$-algebras and has the same residue field  and value group with $\wt{V}$.
  Thus the extension $W\subset \wt{V}$ is immediate. 
  By Theorem \ref{T2} $\wt{V}$ is a filtered union of complete intersection  $W$-algebras and so of complete intersection ${\bf Z}_{(p)}$-algebras as usual using \cite[Lemma 6]{P2} as in Proposition \ref{p}.

Let $E$ be a finitely generated $\bf Z_{(p)}$-algebra and $w:E\to V$ a morphism. Then the composite map $E\to V\to \wt{V}$ factors through a complete intersection ${\bf Z}_{(p)}$-algebra $D$.  Thus  $w$ 
 factors through $D$ too  as in the proof of Proposition \ref{p5}  because all finite systems of polynomial equations over $\bf Z$ which have a solution in $\wt V$ must have one in $V$. This is enough  by \cite[ Lemma 1.5]{S}.
\hfill\ \end{proof}

Let $A\to B\to C$ be two ring morphisms and $H_i(A,B,C)$, $H^i(A,B,C)$ be the Andr\'e-Quillen homology, respectively cohomology (see \cite{An}).

\begin{Proposition} (\cite[X, Corollaire 20]{An}) \label{an} 
Let $(A,\mm,K) $ be a Noetherian local ring. Then $A$ is a complete intersection ring if and only if one of the following conditions holds 

\begin{enumerate}
\item $H_i(A,K,K)=0$ for all $i\geq 3$,
\item $H^i(A,K,K)=0$ for all $i\geq 3$, 
\item $H_3(A,K,K)=0$,
\item $H^3(A,K,K)=0$.
\end{enumerate}
\end{Proposition}

\begin{Proposition} \label{limit} Let $(V,\mm,k)$ be a valuation ring, which is a filtered colimit of complete intersection rings. Then $H_i(V,k,k)=H^i(V,k,k)=0$ for all $i\geq 3$.
\end{Proposition}

\begin{proof} Assume that $V$ is a filtered colimit of some complete intersection local rings $(B_j,\mm_j,k_j)$, $j\in J$. Then 
$H_i(B_j,k_j,k_j)=H^i(B_j,k_j,k_j)=0$ for all $i\geq 3$ and $j\in J$.  Thus 

$H_i(A,k,k)=\varinjlim_{j \in J}H_i(B_j,k_j,k_j)=0$ and

$H^i(A,k,k)=\varinjlim_{j \in J}H^i(B_j,k_j,k_j)=0$

for all $i\geq 3$ using for example \cite[Tag 08S9]{SP}.
\hfill\ \end{proof}

\begin{Corollary} \label{hom} Let $(V,\mm,k)$ be a valuation ring which is in one of the conditions of Theorem \ref{T3}, or
in one of the conditions of Theorem \ref{T4}. Then 
 $H_i(V,k,k)=H^i(V,k,k)=0$ for all $i\geq 3$.
\end{Corollary}

\end{document}